\documentclass[12pt]{article}
\usepackage{amscd,amsfonts,amssymb,amscd,amsmath,latexsym}
\usepackage[hypertex,backref]{hyperref}
\usepackage{pslatex}
\makeatletter
\renewcommand{\subsubsection}{\@startsection
{subsection}
{1}
{0mm}
{0mm}
{0mm}
{\normalfont\normalsize\itshape}}
\makeatother

\newtheorem{theorem}{Theorem}[section]
\newtheorem{prop}[theorem]{Proposition}
\newtheorem{lem}[theorem]{Lemma}
\newtheorem{ddd}[theorem]{Definition}
\newtheorem{kor}[theorem]{Corollary}

\newcommand{\forget}[1]{}

{\catcode`@=11\global\let\c@equation=\c@theorem}

\newcommand{\proof}{{\it Proof.$\:\:\:\:$}}
 
\newcommand{\kaaa}{{\frak k}}
\newcommand{\paaa}{{\frak p}}
\newcommand{\vp}{{\varphi}}
\newcommand{\taaa}{{\frak t}}
\newcommand{\haaa}{{\frak h}}
\newcommand{\uaaa}{{\frak u}}
\newcommand{\R}{\mathbb{R}}

\newcommand{\C}{\mathbb{C}}

\newcommand{\gaaa}{{\frak g}}

\newcommand{\laaa}{{\frak l}}

\newcommand{\maaa}{{\frak m}}
\newcommand{\aaaa}{{\frak a}}
\newcommand{\naaa}{{\frak n}}

\newcommand{\Tr}{{\tt Tr}}

\newcommand{\cZ}{\mathcal{Z}}

\newcommand{\cU}{\mathcal{U}}
\newcommand{\Hom}{{\tt Hom}}

\newcommand{\End}{{\tt End}}

\newcommand{\cF}{\mathcal{F}}
\newcommand{\cD}{\mathcal{D}}

\newcommand{\opp}{{\tt opp}}
\newcommand{\tr}{{\tt tr}}

\newcommand{\Ad}{{\tt  Ad}}

\newcommand{\id}{{\tt id}}

\newcommand{\nat}{\mathbb{N}}
\newcommand{\supp}{{\tt supp}}

\def\imath{{i}}

\def\hB{\hspace*{\fill}$\Box$ \newline\noindent}

\def\hB{\hspace*{\fill}$\Box$ \\[0cm]\noindent}
\newcommand{\cL}{\mathcal{L}}
 
\newcommand{\cug}{{\cU(\gaaa)}}
\newcommand{\Op}{{\tt Op}}

\title{On quantum ergodicity for vector bundles}
\author{Ulrich Bunke and Martin Olbrich
\thanks{Mathematisches Institut, Universit{\"a}t G{\"o}ttingen,
Bunsenstr. 3-5, 37073 G{\"o}ttingen, GERMANY,  bunke@uni-math.gwdg.de,
olbrich@uni-math.gwdg.de}
}

\begin{document}

\maketitle

\begin{abstract}
In the present paper we develop a framework in which questions
of quantum ergodicity for operators acting on sections of hermitian vector bundles over Riemannian manifolds can be studied.
We are particularly interested in the case of locally symmetric spaces.
For locally symmetric spaces, we extend the recent construction of Silberman and Venkatesh \cite{sv04} of representation theoretic lifts to vector bundles.
\end{abstract}

\tableofcontents
\parskip3ex

\section{Introduction}

\subsubsection{}
We start with a brief review of the basic set-up for the study of  quantum ergodicity
of the Laplace operator acting on functions on a Riemannian manifold.
Given a closed Riemannian manifold we can consider a sequence of
normalized eigenfunctions of the Laplace operator associated to a
sequence of eigenvalues tending to infinity. Taking the square of the absolute
value of each of these functions we obtain a sequence of
probability measures on the manifold. 
Note that the space of
probability measures is weakly compact.
So we can ask for a description of possible 
limit points of this sequence. In particular,
in the framework of quantum unique ergodicity,
we want to know under which circumstances 
there is a unique limit point, namely
the measure determined by the Riemannian metric.

\subsubsection{}

The natural way to study these limit measures is to lift them in a canonical
way to a probability measure (called microlocal lift) on the unit sphere bundle of the manifold. One way to define a
microlocal lift is as follows. Each eigenfunction defines a positive
state  on
the algebra of zero-order  pseudodifferential operators.
If we choose a  positivity-preserving operator convention (a right
inverse of the symbol map), then
this state induces a positive linear form on the algebra of
symbols. Since the latter is the algebra of functions on the
unit sphere bundle this linear form is just a measure on this bundle.
Applying this to a sequence of eigenfunctions we get a bounded sequence of
measures on the unit sphere bundle.
Any limit point of this sequence is called a microlocal lift.
It is now an interesting observation that {\em the set of microlocal lifts
is independent of the choice of the operator convention}.

\subsubsection{}\label{plim}

The unit sphere bundle carries a natural dynamical system, the
geodesic flow. It could be considered as the classical counterpart of
the quantum system described by the Laplace operator. The second basic
observation is now that  {\em all microlocal lifts are invariant with
respect to the geodesic flow}. The combination of  this observation
with additional information about
mixing properties of the geodesic flow 
is the starting point of a finer investigation 
of  the shape of these microlocal lifts. In particular,
under the assumption, that the geodesic flow is ergodic (with respect
to Lebesgue class), it is natural to ask wether  the
microlocal lift is just the (normalized) Riemannian measure.
This is the basic question of quantum ergodicity.
We refer to the introduction of \cite{sv04} for a detailed
description of the current knowledge. Here we only mention the following.
A manifold (or rather its Laplacian) for which the Riemannian measure
is the only microlocal lift is called
quantum uniquely ergodic (QUE).
Rudnick and Sarnak \cite{rusar} conjectured that
negatively curved manifolds are always QUE. Recently, Lindenstrauss \cite{lin2}
has proved an arithmetic version of this conjecture for certain arithmetic
hyperbolic surfaces.

\subsubsection{}\label{plum}

The details of the construction of microlocal lifts and the verification
of the two basic properties are not at all complicated. It is the
purpose of Section \ref{gen} to give these arguments in a more general
setting. In fact, if the Riemannian manifold comes
equipped
with a hermitian vector bundle with connection, then we can replace
the Laplace operator on the manifold by the Laplace operator on this
bundle. Then we are looking for microlocal lifts associated with
sequences of eigensections of the operator.
The new point is that the algebra of symbols is now the
algebra of sections of the endomorphism bundle of the vector bundle
lifted to the unit sphere bundle. In particular, this algebra can be  non-commutative.
This essentially leads to a change of terminology, the main instance
of which is the replacement of probability measures by states.

The set of microlocal lifts is now a set of states on the algebra of
symbols. We show in Proposition \ref{vvd} that this set is naturally associated to the
geometric data. The connection induces a natural lift of the geodesic
flow to a flow of automorphisms of the algebra of symbols, and we
verify in Proposition \ref{inv1} that each microlocal lift is invariant.

Finer quantum ergodicity questions are left untouched in this paper
and will be a topic of future research.
Note that in the bundle case one cannot expect the microlocal lift
to be  unique (even for
negatively curved manifolds). 
Let us consider e.g. the case of differential forms which can 
be decomposed into closed and coclosed ones. 
Associated with this decomposition is 
a
natural splitting 
into two parts of the pull-back of bundle of
differential forms to the unit sphere
bundle.
If we consider e.g  sequence of closed eigenforms, then
the associated microlocal lifts are annihilated by the projection onto
the  subbundle corresponding to coclosed eigenforms.
The microlocal lifts associated to sequences of coclosed forms behave
in the opposite way.

\subsubsection{}

In  Section \ref{gf} we start to develop a theory of representation
theoretic lifts for the case of 
a compact locally symmetric space $\Gamma\backslash G/K$.
Representation theoretic lifts serve as a substitute for the microlocal lifts
discussed so far. They are designed  to take  into account the rich structure available in the locally symmetric situation. While defined without any reference to pseudodifferential
operators it turns out (see the final Section \ref{luzio}) that they determine the microlocal lifts. Thus representation theoretic lifts should be considered
as refined microlocal ones.

Our guide here is the recent paper by Silberman and Venkatesh
\cite{sv04}, where the notion of a
representation theoretic lift was introduced.
The Laplace operator on functions on $\Gamma\backslash G/K$ commutes with a whole algebra of differential operators coming from
the center $\cZ(\gaaa)$  of the universal enveloping algebra of the Lie algebra
of $G$. Therefore the spectrum of the Laplace operator can be
further decomposed with respect to this algebra. For studying the fine
structure of the lifts associated with the locally symmetric situation it seems more
appropriate to consider instead of a sequence of eigenfunctions of the
Laplacian the corresponding sequence of embeddings of spherical unitary representations of $G$  into $L^2(\Gamma\backslash G)$.
Following earlier constructions for special cases due to Zelditch and Lindenstrauss, Silberman and Venkatesh associate to such a sequence of embeddings a representation theoretic lift.

Our main observation is that one can apply an analogous procedure if one
wants to study sequences of eigensections of bundles of the form
$\Gamma\backslash G\times_K V_\gamma$, where $(\gamma,V_\gamma)$ is
a unitary representation of $K$.

\subsubsection{}
Let us remark at this point that locally symmetric spaces of higher rank do not have strictly negative
curvature. In fact, they do not have the QUE-property
defined in \ref{plim} as follows from
the results of Section \ref{luzio}. 
In fact it turns out that the microlocal lifts associated to
conveniently arranged sequences
(see \ref{cas})
of embeddings of principal series representations
are supported on subsets
of the
unit sphere bundle of Lebesgue measure zero. 

Thus the definition of QUE must be modified in the case of higher rank
locally symmetric spaces. Replacing microlocal by representation theoretic lifts which live
on $\Gamma\backslash G$ one can define an (arithmetic)
QUE-property which has  recently  been verified 
in many cases (see the forthcoming second part of \cite{sv04}).

\subsubsection{}

We now describe our construction  of the representation theoretic
lifts which is presented in detail in Section \ref{gf}.
We follow quite closely the approach of \cite{sv04}. Our
  contribution is essentially an adaptation of arguments and
  language to the
  non-commutative situation in the case of non-trivial $K$-types.

The main part of the spectrum of the Laplacian (and the other locally invariant
differential operators) on $L^2(\Gamma\backslash G\times_K V_\gamma)$ is caused by
embeddings of unitary principal series representations
(associated with the minimal parabolic
subgroup of $G$) into $L^2(\Gamma\backslash G)$. These pricipal series representations
come in natural families.
Let $G=KAN$ be an Iwasawa decomposition of $G$ and $M:=Z_K(A)$ be the
 centralizer of $A$ in $K$. Then a family of unitary principal series
 representations (see \ref{hs}) is determined by an element $\kappa\in \hat M$.
The family parameter runs through $\aaaa^*$, the dual of the Lie algebra
of $A$.

\subsubsection{}

We fix now $\kappa$ and therefore a  family of principal series
representations. We consider
 a sequence of embeddings of members of this family into $L^2(\Gamma\backslash
G)$  
 with parameter tending to infinity approximately along a regular ray
 (see \ref{reg})
 in $\aaaa^*$.
Such a sequence contributes to the spectrum in $L^2(\Gamma\backslash G\times_K V_\gamma)$ if and only if $[V_\kappa\otimes V_\gamma]^M\ne\{0\}$.
In the
function case ($\gamma=1$) this condition is equivalent to $\kappa=1$, i.e.,
the corresponding principal series representations are spherical.
In contrast
to the spherical case the dimension of $[V_\kappa\otimes V_\gamma]^M$
can be greater  than one in general.

To any $T\in
 [V_\kappa\otimes V_\gamma]^M$ there is  an associated  vector $\psi_T$ of ``$\gamma$-spherical elements''
in the corresponding principal series representation (see \ref{hs}).
 Applying the embedding $\xi$ of a principal series
representation  into $L^2(\Gamma\backslash
G)$ to $\psi_T$ we get an eigensection $\xi(\psi_T)\in L^2(\Gamma\backslash G\times_K V_\gamma)$ denoted by $\xi(\psi_T)$.
This section defines a state $\sigma_{\xi(\psi_T)}$  on the algebra
$C(\Gamma\backslash G\times_K \End(V_\gamma))$,
which can be identified
with the subalgebra of $K$-invariants in $C(\Gamma\backslash G)\otimes
\End(V_\gamma)$. Recall that we considered a sequence of  embedded
principal series representations, and therefore we get a sequence of
such states. We are interested in ``lifting'' the limit states to
a state of $C(\Gamma\backslash G)\otimes
\End(V_\gamma)$ which is not a priori $K$-invariant.

To this end we construct for each individual embedding $\xi$ a functional
 $\sigma^{\xi}_{\psi_T,\delta_T}$ (see \ref{spaet} for the definiton) on the
 smaller algebra
 $C^K(\Gamma\backslash G)\otimes \End(V_\gamma)$ (the $C^K$ stands for
 $K$-finite functions). After choosing an appropriate subsequence
of embeddings the corresponding sequence of functionals converges to a functional which extends to
a state on the algebra $C(\Gamma\backslash G)\otimes \End(V_\gamma)$ (see
Proposition \ref{pos}).
The states which are obtained in this way are the representation
theoretic lifts in question. Their restriction to the subalgebra of $K$-invariants in $C(\Gamma\backslash G)\otimes
\End(V_\gamma)$ coincides with the set of limit states associated to
$\sigma_{\xi(\psi_T)}$ discussed above (see Proposition
\ref{fr}). On the one hand, this 
justifies the name ``lift''. On the other hand, the representation
theoretic lifts contain additional microlocal information.

\subsubsection{}

It is not apriori clear that a representation theoretic lift
is a limit point of states associated to a sequence of functions in
$L^2(\Gamma\backslash G)\otimes V_\gamma$. In Theorem \ref{dde} we show
that this is indeed the case.

The space of embeddings of a fixed unitary
representation of $G$ into $L^2(\Gamma\backslash G)$  is acted on by Hecke operators.
In particular, fixing a Hecke operator and a complex number,
it makes sense to talk about eigenembeddings with this given eigenvalue.

If the representation theoretic lift
is associated to a family of eigenembeddings (for fixed Hecke operator
and eigenvalue), then we show further that the representation theoretic
lift is a limit of states associated to a  sequence of functions in
$L^2(\Gamma\backslash G)\otimes V_\gamma$ which are also eigenfunctions
of the Hecke operator with the given eigenvalue.
Such a property played a fundamental role in the study of the
quantum ergodic properties
in \cite{lin}. In particular, it implies restrictions on the support of
the representation theoretic lifts.

\subsubsection{}

In Section \ref{dgl} we investigate various invariance properties
of the representation theoretic lifts.
Note that the group $A$ acts on $\Gamma\backslash G$ by
right multiplication.
This induces an action of $A$ by automorphisms on
$C(\Gamma\backslash G)\otimes
\End(V_\gamma)$.
The main result of the section (Theorem \ref{inv}) now states that
{\em all representation theoretic lifts are invariant with respect to
  the action of $A$}.
This is the locally symmetric counterpart of the invariance of the
microlocal lifts with respect to the geodesic flow.

We show further
(Proposition \ref{mulm}) that the microlocal lifts are $M$-invariant and,
when associated with $T\in
 [V_\kappa\otimes V_\gamma]^M$, are essentially
states on the smaller algebra $C(\Gamma\backslash G\times_M \End(V_T))$,
where the irreducible $M$-representation space $V_T\subset V_\gamma$
is given by $\{<v,T>\:|\: v\in V_\kappa\}$. In particular, a pair of sequences
of eigensections which corresponds to a pair of linearly independent
$T$'s has a disjoint pair of sets of representation theoretic lifts. In view of the comparison
result (Corollary \ref{rg3}) between representation theoretic and microlocal lifts obtained
in the final Section \ref{luzio} this is another manifestation of the non-uniqueness
of the microlocal lifts mentioned at the end of \ref{plum}.

\section{Microlocal lifts}\label{gen}

\subsubsection{}

Let $M$ be a closed smooth manifold with Riemannian metric $g$. Let
$E\rightarrow M$ be a complex vector bundle of dimension $n$ with
hermitian metric $h$ and metric connection $\nabla$. By
$\Delta=\nabla^*\nabla$ we denote the Laplace operator.
With the smooth sections $C^\infty(M,E)$  as domain it can be
considered as an unbounded essentially selfadjoint operator on the
Hilbert space $L^2(M,E)$ of square integrable sections of $E$.

\subsubsection{}

 Let $\pi:SM\rightarrow M$ denote the unit sphere bundle of the cotangent
bundle $T^*M\rightarrow M$.
Let $\Psi D O^i(M,E)$ denote the algebra of classicial $i$'th-order
pseudodifferential operators on $M$.
Then we have an exact sequence of algebras
$$0\rightarrow \Psi D O^{-1}(M,E)\rightarrow \Psi D
O^0(M,E)\stackrel{s}{\rightarrow}
C^\infty(SM,p^*\End(E))\rightarrow 0\ ,$$
where $s$ is the principal symbol map.

A linear continuous right-inverse $\Op:C^\infty(SM,p^*\End(E))\rightarrow \Psi D
O^0(M,E)$  of $s$ is called a quantization or
operator convention. The algebra pseudo-differential operators is here
topologized as algebra of continuous operators on $C^\infty(M,E)$.
In general, a quantization does not extend continuously to a map 
$C(SM,p^*\End(E))\rightarrow B(L^2(M,E))$. But using a
construction of Friedrichs (see \cite{taylor81}, p. 142) we can choose the
quantization $\Op$ such that it preserves positivity, i.e.
if $a\in C^\infty(SM,p^*\End(E))$ is a non-negative element in the $C^*$-algebra $C(SM,p^*\End(E))$,
then $\Op(a)\ge 0$ in the $C^*$-Algebra $B(L^2(M,E))$.

\subsubsection{}

Consider a $C^*$-algebra $A$ and a dense subalgebra $A_\infty\subset A$.
A linear map $\sigma:A_\infty\rightarrow \C$ is called positive if
$a\ge 0$ implies that $\sigma(a)\ge 0$.

If $\sigma:A_\infty\rightarrow \C$ is positive and $\sigma(1)<\infty$, then $\sigma$ extends
uniquely to a continuous linear positive map $\sigma:A\rightarrow
\C$. A state on $A$ is a normalized (i.e. $\sigma(1)=1$) linear positive map $\sigma:A\rightarrow
\C$.

%A weight $\sigma$ on a $C^*$-algebra $A$ is a function
%$\sigma:A^+\rightarrow[0,\infty]$ which is positive homogeneous
%and additive, where $A^+\subset A$ is the cone of positive elements of
%$A$. If $\{a\in A^+|\sigma(a)<\infty\}\subset A$ is dense, then
% $\sigma$ is  densely defined.
%A weight is normalized (bounded) if $\sigma(1)=1$
%($\sigma(1)<\infty$). A densely defined normalized (bounded)
%weight extends 
%to a state (continuous positive linear functional) $\sigma:A\rightarrow \C$.  

\subsubsection{}\label{statdef}

Let $\psi\in C^\infty(M,E)$ be a unit vector in $L^2(M,E)$.
We then consider the linear map
$$\sigma_\psi:C^\infty(SM,p^*\End(E))\rightarrow \C$$
given by 
$$\sigma_\psi(a)=<\psi,\Op(a)\psi>\ .$$
Since the quantization preserves positivity,
$\sigma_\psi$ is a positive.
Since $\sigma_\psi(1)<\infty$ it follows that
$\sigma_\psi$  extends to a continuous
positive linear functional
$$\sigma_\psi:C(SM,p^*\End(E))\rightarrow \C\ .$$
In fact, we have the uniform estimate
$$\|\sigma_\psi\|\le \|\Op(1)\|\ .$$

\subsubsection{}

We now consider the countable set of functionals
on $C(SM,p^*\End(E))$ of the form $\sigma_\psi$,
where $\psi$ is a normalized eigenvector of $\Delta$.
This set is a bounded set in the Banach dual of $C(SM,p^*\End(E))$
and therefore weak-$*$-precompact. By $V\subset C(SM,p^*\End(E))^*$ we
denote then non-empty set of all its accumulation points.

\begin{prop}\label{vvd}
The set $V$ is independent of the choice of the positive quantization map and
consists of states.
\end{prop}
\proof
Let $V^\prime$ denote the set defined with another choice $\Op^\prime$
of the quantization. For the first assertion it suffices to show that $V\subset V^\prime$.
We consider $\sigma\in V$. Then there exists a sequence of normalized eigenvectors
$\psi_n$ of $\Delta$ to eigenvalues $\lambda_n\to\infty$ such that  
 for all $f\in C^\infty(SM,p^*\End(E))$ we have
$\sigma(f)=\lim <\psi_n,\Op(f)\psi_n>$.
 Note that $\Op(f)-\Op^\prime(f)\in \Psi DO^{-1}(M,E)$.
But for $A\in \Psi DO^{-1}(M,E)$ we have
$\lim <\psi_n,A\psi_n>=0$ since $\lambda_n\rightarrow \infty$.
This shows that $$\sigma(f)=\lim <\psi_n,\Op^\prime(f)\psi_n>\ .$$
We conclude that $\sigma\in V^\prime$.

We now show that $V$ consists of states.
It is clear that the elements of $V$ are positive. We must verify normalization.
Note that $\Op(1)=1+A$, where $A\in\Psi DO^{-1}(M,E)$.
We conclude that $\sigma(1)=1$. \hB

\subsubsection{}

Using functional calculus we can define the strongly continuous group
of unitary operators $\exp(it\sqrt{\Delta})$. One can show that these
are Fourier-integral operators. Note that the conjugation by a Fourier-integral
operator preserves pseudodifferential operators.
Let $f\in C^\infty(SM,p^*\End(E))$. Then
in principle one can calculate the symbol
of $$\exp(it\sqrt{\Delta})\Op(f)\exp(-it\sqrt{\Delta})$$ using the
calculus of Fourier-integral operators. Here we prefer a simpler way by computing the infinitesimal
action which amounts to a calculation of the symbol
$s([i\sqrt{\Delta},\Op(f)])\in C^\infty(SM,p^*\End(E))$.
Let $X\in C^\infty(SM,TSM)$ denote the generator of the geodesic flow.
Note that we have an induced connection on $p^*\End(E)$ which we will
also denote by $\nabla$.

\begin{lem}\label{infac}
We have $s([i\sqrt{\Delta},\Op(f)])=\nabla_X f$.
\end{lem}
\proof
This is a local computation and independent of the choice of the
quantization.

We consider a point in $M$ and choose geodesic normal coordinates
$x$. Let $(x,\xi)$ denote the corresponding coordinates of $T^*M$.
We want to compute $s([i\sqrt{\Delta},\Op(f)])$ in the point $(0,\xi)$ with
$\|\xi\|=1$. We further trivialize $E$ using radial parallel
transport. We now use the standard quantization.
In these coordinates the full symbol of $\sqrt{\Delta}$ is given by
$\|\xi\|+O(x^2)$. The principal symbol of a commutator of pseudodifferential
operators is given by the Poisson bracket of their symbols. 
Therefore we get (on the sphere $\{\|\xi\|=1\}$)
$$s([i\sqrt{\Delta},\Op(f)])(0,\xi)=\xi^i \partial_{x_i} f(0,\xi)\ .$$
This implies the assertion in view of the choice of the
trivializations,
since $\xi^i\partial_{x_i}$ is the value of the generator of the
geodesic flow at $(0,\xi)$.
\hB

\subsubsection{}

Using the connection $\nabla$ on $p^*\End(E)$ we can lift the geodesic
flow $\Phi_t$ on $SM$ to a flow $\tilde \Phi_t$ on $p^*\End(E)$. 
We denote the action of this flow on sections by the same symbol. We
have $\frac{d}{dt}_{|t=0}\tilde \Phi_t(f)=\nabla_Xf$. 
Thus the algebra $C(SM,p^*\End(E))$ comes with a flow of
automorphisms $\tilde \Phi_t$. By Lemma \ref{infac} we have
$$s(\exp(it\sqrt{\Delta})\Op(f)\exp(-it\sqrt{\Delta}))=\tilde \Phi_{t}(f)\ .$$

\begin{prop}\label{inv1} 
Every limit state $\sigma\in V$ is invariant under this flow.
\end{prop}
\proof
Let $\psi_n$ be a sequence of normalized eigenvectors
of $\Delta$ to eigenvalues $\lambda_n$ such that
for all $f\in C^\infty(SM,p^*\End(E))$ we have
$\sigma(f)=\lim <\psi_n,\Op(f)\psi_n>$.
We compute 
\begin{eqnarray*}
\sigma(f)&=&\lim <\psi_n,\Op(f)\psi_n>\\
&=&\lim <\exp(-it\sqrt{\Delta})\psi_n,,\Op(f)\exp(-it\sqrt{\Delta})\psi_n>\\
&=&\lim <\psi_n,\exp(it\sqrt{\Delta})\Op(f)\exp(-it\sqrt{\Delta})\psi_n>\\
&=&\lim <\psi_n,\Op(s(\exp(it\sqrt{\Delta})\Op(f)\exp(-it\sqrt{\Delta})))\psi_n>\\
&=&\sigma(\tilde \Phi_{t}(f))\ .
\end{eqnarray*}
\hB

\subsubsection{}

A state on the algebra of functions $C(SM)$ is the same
thing as
a probability measure on $SM$. A state $\sigma$ on the algebra
$C(SM,p^*\End(E))$
determines and is  determined by a pair $(\mu,M)$, where
$\mu$ is a probability measure on $SM$, $M\in
L^\infty(SM,p^*\End(E),\mu)$  gives
a measurable family of states on the local algebras, and
$\sigma(f)=\mu(\tr Mf)$.
Here $\tr:p^*\End(E)\rightarrow \C$ denotes the local trace and
$\tr Mf\in L^\infty(SM,\mu)$.

If $\sigma$ is invariant under the flow $\tilde \Phi_t$, then $\mu$ is invariant under
the geodesic flow, and $M$ is invariant under its lift $\tilde \Phi_t$.

\subsubsection{}

The picture which we have described so far is a simple generalization
of a well-known construction (see \cite{scn}, \cite{cv}, \cite{zel})
from the case of the trivial bundle $E=M\times \C$ to arbitrary
bundles $E\rightarrow M$.
It is by now an interesting piece of mathematics to obtain more
information
about the size of the set of limiting states $V$ and the properties
of its elements under ergodicity assumptions on the geodesic flow
$\Phi_t$.

\section{Representation theoretic lifts}\label{gf}

\subsubsection{}

Let $G$ be a semisimple Lie group, $K\subset G$ be a maximal compact
subgroup of $G$ and $\Gamma$ be a cocompact torsion free discrete
subgroup. Then we can consider the locally symmetric space $M=\Gamma\backslash
G/K$ with a Riemannian metric given by the Killing form of $G$.  A
unitary representation $(\gamma,V_\gamma)$ of $K$ gives rise to a vector
bundle $V(\gamma):=\Gamma\backslash G\times_K V_\gamma$ over $M$ which comes
with a natural connection. In this situation we can consider the limit
states as discussed in Subsection \ref{gen}. But because of the
locally-symmetric structure we can perform a refined construction which
we will describe below.

\subsubsection{}

We consider the $C^*$-algebra of functions $A_\gamma:=C(\Gamma\backslash
G)\otimes \End(V_\gamma)$. Let $(\pi,H)$ be an irreducible unitary representation
of $G$. By $H_{\pm\infty}$ we 
denote the distribution or smooth vectors of $H$. Let
 $\xi\in [H^*_{-\infty}]^\Gamma$ be an invariant
distribution vector. Equivalently, we can consider $\xi$ as an
embedding
$H_\infty\rightarrow C^\infty(\Gamma\backslash G)$. 
We shall assume that $\xi$ is normalized such that this embedding
extends to a unitary embedding $H\rightarrow L^2(\Gamma\backslash G)$. 
Let $\phi,\psi\in H_\infty\otimes V_\gamma$. We have $\xi(\phi)\in
C^\infty(\Gamma\backslash G)\otimes V_\gamma$.
Then we can define a functional
$\sigma_{\phi,\psi}$ on $A_\gamma$ by
\begin{equation}\label{sswe}\sigma^\xi_{\phi,\psi}(f)=\int_{\Gamma\backslash G}
<\xi(\phi)(\Gamma g),f(\Gamma g)\xi(\psi)(\Gamma g)>\ .\end{equation}
If $\|\phi\|=1$, then $\sigma^\xi_{\phi,\phi}$ is  a  state.

\subsubsection{}

Let $\phi,\psi\in H_\infty$ and $f\in A_\gamma^\infty:=C^\infty(\Gamma\backslash
G)\otimes \End(V_\gamma)$. Furthermore, let $\gaaa$ denote the Lie
algebra of $G$ and $X\in \gaaa$.
We set $Xf(g)=f(gX):=\frac{d}{dt}_{|t=0}f(g\ \exp(tX))$.
\begin{lem}\label{leib}
We have
$$\sigma^\xi_{\pi(X)\phi,\psi}(f)+\sigma^\xi_{\phi,\pi(X)\psi}+\sigma^\xi_{\phi,\psi}(Xf)=0\
.$$
\end{lem}
\proof
This follows from the fact that the $G$-action on $\Gamma\backslash G$
preserves the measure. \hB

\subsubsection{}\label{hs}

Let $G=KAN$, $g=k(g)a(g)n(g)$, be an Iwasawa decomposition, and let $M:=Z_K(A)\subset K$
be the centralizer of $A$ in $K$. Let $(\kappa,V_\kappa)$ be an
irreducible unitary representation of $M$.

Let $\aaaa$ denote the Lie algebra of $A$.
The choice of $\kappa$ gives rise to a family of unitary principal series
representations $(\pi^{\kappa,i\lambda},H^{\kappa,i\lambda})$,
$\lambda\in \aaaa^*$, of $G$. In the compact picture we set
$H^{\kappa,i\lambda}:=L^2(K\times_MV_\kappa)$.
Then $(\pi^{\kappa,i\lambda}(g)\phi)(k)=\phi(k(g^{-1}k))
a(g^{-1}k)^{i\lambda-\rho}$. Here  for $a\in A$ and $\lambda\in
\aaaa_\C:=\aaaa\otimes_\R\C$ the symbol $a^\lambda$ is a short-hand for
$\exp(\lambda(\log(a)))$.
Moreover, $\rho\in \aaaa^*$ is given by
$\rho(H)=\frac{1}{2}\Tr\, \Ad(H)_{|\naaa}$, $H\in \aaaa$, where $\naaa$
denotes the Lie algebra of $N$.

By Frobenius reciprocity we have
$[C^\infty(K\times_MV_\kappa)\otimes V_\gamma]^K\cong [V_\kappa\otimes
V_\gamma]^M$. Explicitly the isomorphism is given by evaluation at
$1\in K$. For $T\in [V_\kappa\otimes
V_\gamma]^M$ we let $\psi_T\in [C^\infty(K\times_MV_\kappa)\otimes V_\gamma]^K$ denote the corresponding
element.

%If $\xi:H^{\kappa,i\lambda}\rightarrow L^2(\Gamma\backslash G)$ is an
%unitary embedding and $T,T^\prime\in  [V_\kappa\otimes
%V_\Gamma]^M$, then we define the state
%$$\sigma^\xi_{T,T^\prime}:=\sigma^\xi_{\psi_T^\lambda,\psi_{T^\prime}^\lambda}\
%$$ on $A$.

\subsubsection{}\label{reg}

The Weyl group $W(G,A):=N_K(A)/M$ acts by reflections on $\aaaa^*$.
A point $\lambda\in\aaaa^*$ is called singular, if it is fixed by some
element of the Weyl group. Otherwise it is called regular.
Let $S(\aaaa^*):=(\aaaa^*\setminus\{0\})/\R^+$ be the space of rays in $\aaaa^*$.
We have a corresponding decomposition of $S(\aaaa^*)$ into singular
and regular rays. Using the metric $\|.\|$ induced by the Killing form
on $\gaaa$ we will identify $S(\aaaa^*)$ with the unit sphere in
$\aaaa^*$.

\subsubsection{}

For $0\not=\lambda\in \aaaa^*$ let $[\lambda]\in S(\aaaa^*)$ denote
the corresponding ray.

\begin{ddd}
We define the closed subset $L(\kappa)\subset S(\aaaa^*)$ as the set
points $l\in S(\aaaa^*)$ such that there exists
a sequence $\xi_n:H^{\kappa,i\lambda_n}\rightarrow
L^2(\Gamma\backslash G)$  of 
unitary embeddings such that $\lambda_n\rightarrow \infty$ and
$[\lambda_n]\rightarrow l$.
\end{ddd}

\subsubsection{}

Let $C^K(K\times_M V_\kappa)\subset C^\infty(K\times_M V_\kappa)$
denote the subspace of $K$-finite vectors.
Let us fix a regular point $l\in L(\kappa)$.
\begin{ddd}\label{cas} We call a sequence $\xi_n:H^{\kappa,i\lambda_n}\rightarrow
L^2(\Gamma\backslash G)$  of 
unitary embeddings $l$-conveniently arranged if $\lambda_n\rightarrow \infty$,
$[\lambda_n]\rightarrow l$, all $\lambda_n$ are regular,
and for all
$\phi,\psi\in C^K(K\times_M V_\kappa)\otimes V_\gamma$
the sequence of functionals $\sigma^{\xi_n}_{\phi,\psi}$
converges weakly.
\end{ddd}

\subsubsection{}

\begin{lem}
For each regular $l\in L(\kappa)$ there exists a $l$-conveniently arranged
sequence.
\end{lem}
\proof
Consider a sequence $\xi_n:H^{\kappa,i\lambda_n}\rightarrow
L^2(\Gamma\backslash G)$  of 
unitary embeddings such that $\lambda_n\rightarrow \infty$ and
$[\lambda_n]\rightarrow l$. 
By taking a subsequence we can assume that all $\lambda_n$ are
regular.

For fixed $\phi,\psi\in  C^K(K\times_M V_\kappa)\otimes V_\gamma$
the bounded set $\{\sigma^{\xi_n}_{\phi,\psi}|n\in\nat\}$ of functionals on $A_\gamma$ is
weak-$*$-precompact. Therefore, by taking a subsequence, we can assume
that $\sigma^{\xi_n}_{\phi,\psi}$ weakly converges.

Finally note that $C^K(K\times_M V_\kappa)\otimes V_\gamma$ is a
vector space with a countable base. Therefore by a diagonal sequence argument
we can again choose a subsequence such that
 $\sigma^{\xi_n}_{\phi,\psi}$ converges
 for all $\phi,\psi\in  C^K(K\times_M V_\kappa)\otimes V_\gamma$.
\hB

%For each point $l\in L(\kappa)$ and $\psi,\phi\in  C^\infty(K\times_MV_\kappa)\otimes V_\gamma$ we define the following set
%of states $V(\psi,\phi,l)$ on $A_\gamma$.

%\begin{ddd}
%A state $\sigma$ on $A$ belongs to $V(\phi,\psi,l)$ if there
%exists a sequence $\xi^n:H^{\kappa,i\lambda_n}\rightarrow
%L^2(\Gamma\backslash G)$ of 
%unitary embeddings such that $\lambda_n\rightarrow \infty$ and
%$[\lambda_n]\rightarrow l$, and such that
%$\sigma$ is the weak limit of the sequence $\sigma^{\xi_n}_{\phi,\psi}$.
%\end{ddd}

\subsubsection{}

Note that the $C^*$-algebra $C(K/M)$ acts naturally on the Hilbert
space $L^2(K\times_M V_\kappa)\otimes V_\gamma$ by multiplication
on the first factor. The subalgebra $C^K(K/M)$ of $K$-finite functions
acts on $C^K(K\times_M V_\kappa)\otimes V_\gamma$. 

We consider a regular $l\in L(\kappa)$ and let $\xi_n$ be a $l$-conveniently arranged sequence as in Definition
\ref{cas}. By $\sigma_{\phi,\psi}$ we denote the weak limit of the
sequence of functionals
$\sigma^{\xi_n}_{\phi,\psi}$ for fixed $\phi,\psi\in C^K(K\times_M
V_\kappa)\otimes V_\gamma$. 

Note that if  $\|\phi\|=1$, then $\sigma_{\phi,\phi}$ is a state.

\subsubsection{}

Let $h\in C^K_\R(K/M)$ be a real-valued
$K$-finite function.
\begin{lem}\label{45d}
We have $$\sigma_{h\phi,\psi}=\sigma_{\phi,h\psi}\ .$$
\end{lem}
\proof
Let $P=MAN\subset G$ be the minimal parabolic subgroup of $G$.
For $\lambda\in\aaaa_\C^*$ we consider the representation
$\kappa_\lambda(man):=\kappa(m)a^{\rho-\lambda}$ of $P$ on $V_\kappa$.
We identify $C^\infty(K\times_MV_\kappa)\otimes V_\gamma$ with
$C^\infty(G\times_P V_{\kappa_\lambda})\otimes V_\gamma$ by
restriction from $G$ to $K$. This restriction intertwines the
$G$-action $\pi^L$ on $C^\infty(G\times_P V_{\kappa_\lambda})\otimes V_\gamma$
by left translations with the action $\pi^{\kappa,\lambda}$ (see
\ref{hs}) on $C^\infty(K\times_MV_\kappa)\otimes V_\gamma$.

Let $\gaaa=\kaaa\oplus\aaaa\oplus\naaa$ be the Iwasawa decomposition,
and let $X=X_\kaaa+X_\aaaa+X_\naaa$ denote the corresponding
decomposition of $X\in \gaaa$.
Furthermore, for $k\in K$ let $X(k):=\Ad(k^{-1})X$ and Iwasawa
decompose $X(k)=X_\kaaa(k)+X_\aaaa(k)+X_\naaa(k)$. 
If $\phi\in C^\infty(G\times_P V_{\kappa_\lambda})\otimes V_\gamma$, then we have
\begin{eqnarray*}
(\pi^{L}(X)\phi)(k)&=&
\phi(Xk)\\&=&\phi(kX(k))\\
&=&\phi(kX_\kaaa(k))+\phi(kX_\aaaa(k))+\phi(kX_\naaa(k))\\
&=& \phi(kX_\kaaa(k))-(\rho-\lambda)(X_\aaaa(k))\phi(k) 
\end{eqnarray*}

Let now $f\in A_\gamma^\infty$ and $\phi,\psi\in C^K(K\times_M
V_\kappa)\otimes V_\gamma$. Note that the action of $\gaaa$ preserves $K$-finite
vectors. Since
$K\ni k\mapsto (\rho-\lambda)(X_\aaaa(k))=:p_{X,\rho-\lambda}(k)\in \C$
is a $K$-finite function on $K/M$, we see that $k\mapsto \phi_X(k):=\phi(kX_\kaaa(k))$
is $K$-finite, too.

By Lemma \ref{leib}
we have
$$0=\sigma^{\xi_n}_{\pi^{L}(X)\phi,\psi}(f)+\sigma^{\xi_n}_{\phi,\pi^{L}(X)\psi}(f)+\sigma^{\xi_n}_{\phi,\psi}(Xf)=0\
.$$
This implies
\begin{eqnarray*}
\sigma^{\xi_n}_{p_{X,i\frac{\lambda_n}{\|\lambda_n\|}}\phi,\psi}(f)+
\sigma^{\xi_n}_{\phi,p_{X,i\frac{\lambda_n}{\|\lambda_n\|}}\psi}(f)&=&
\frac{1}{\|\lambda_n\|}\big(\sigma^{\xi_n}_{p_{X,\rho}\phi,\psi}(f)+
\sigma^{\xi_n}_{\phi,p_{X,\rho}\psi}(f)\\&&-\sigma^{\xi_n}_{\phi_X,\psi}(f)-
\sigma^{\xi_n}_{\phi,\psi_X}(f)-\sigma^{\xi_n}_{\phi,\psi}(Xf)\big)\ .
\end{eqnarray*}
The right-hand side of this equation converges to zero as $n$ tends to
infinity. We set  $$p_{X,l}(k):=-i\lim_{n\to\infty}
p_{X,i\frac{\lambda_n}{\|\lambda_n\|}}(k)\ .$$
The sequences of functions $p_{X,i\frac{\lambda_n}{\|\lambda_n\|}}\phi$ and
$p_{X,i\frac{\lambda_n}{\|\lambda_n\|}}\psi$ span 
finite-dimensional spaces. Since $\sigma_{\phi,\psi}$ is conjugated
linear in $\phi$ 
we conclude that 
$$\sigma_{p_{X,l}\phi,\psi}(f)-\sigma_{\phi,p_{X,l}\psi}(f)=0\
    .$$
 
Let $\cF\subset C^K_\R(K/M)$ denote the algebra of functions generated by
the constant functions and the functions $p_{X,l}$, $X\in \gaaa$.
Then we have shown that for all $h\in \cF$ we have
$$\sigma_{h\phi,\psi}(f)=\sigma_{\phi,h\psi}(f)\
    .$$

It remains to show that $\cF=C_\R^K(K/M)$.
It suffices to show that $\bar \cF=C(K/M)$, where $\bar \cF$ is the
closure of $\cF$ in $C_\R(K/M)$.
The algebra $\cF$  
contains the identity and separates points.
This can be seen as follows.
Using the Iwasawa decomposition we extend $l$ to $\tilde l\in
\gaaa^*$. Then we can write $p_{X,l}(k)=\tilde
l(\Ad(k^{-1})X)=\Ad(k)(\tilde l)(X)$. If $p_{X,l}(k_1)=p_{X,l}(k_2)$
for all $X\in\gaaa$, then we have 
$\Ad(k_1)(\tilde l)=\Ad(k_2)(\tilde l)$.
Since $l$ is regular, this implies that $k_1^{-1}k_2\in M$, hence
$k_1M=k_2M$. We conclude by the Stone-Weierstrass theorem that
$\bar \cF=C_\R(K/M)$.
\hB

\begin{kor}\label{ds}
For 
$\phi,\psi\in C^K(K\times_M
V_\kappa)\otimes V_\gamma$ and $h\in C^K(K/M)$ we have
$$\sigma_{h\phi,\psi}=\sigma_{\phi,\bar h\psi}\ .$$
\end{kor}

\subsubsection{}\label{dif}

Let $A_{\gamma,K}\subset A_\gamma$ denote the subalgebra of $K$-finite
elements. Let $\phi\in C^K(K\times_MV_\kappa)\otimes
V_\gamma$ and $\psi\in C^{-\infty}(K\times_MV_\kappa)\otimes
V_\gamma$. Then we have
$\psi=\sum_{\mu\in \hat K} \psi_\mu$ in the sense of distributions, where
$\psi_\mu\in C^K(K\times_MV_\kappa)\otimes
V_\gamma$ is the component of $\psi$ in the $\mu$-isotypic
subspace of $C^{-\infty}(K\times_MV_\kappa)$.
For $F\subset \hat K$ we set $\psi_F:=\sum_{\mu\in F} \psi_\mu$.
 It is easy to see that
for $f\in A_{\gamma,K}$ the sum
$\sum_{\mu\in \hat K} \sigma_{\phi,\psi_\mu}(f)$
is finite. In fact, there exists a finite set $F=F(\phi,f)$, 
independent of $\psi$ (but which depends on $\phi$ and $f$), of $K$-types which can contribute to this sum.
\begin{ddd}
We define the functional
$\sigma_{\phi,\psi}:A_{\gamma,K}\rightarrow \C$ by
$$\sigma_{\phi,\psi}(f):=\sum_{\mu\in \hat K}
\sigma_{\phi,\psi_\mu}(f)\ .$$
\end{ddd}

\subsubsection{}\label{sit}

Now we fix $T\in [V_\kappa\otimes V_\gamma]^M$ such that $\|T\|=1$ and let $\psi_T\in
C^K(K\times_MV_\kappa)\otimes V_\gamma$ as in \ref{hs}.
Furthermore, we let $\delta_T\in C^{-\infty}(K\times_MV_\kappa)\otimes
V_\gamma$ be the distribution
$\phi\mapsto <T,\phi(1)>$.

\begin{ddd}\label{de1}
We define
 the functional $\sigma_T$ on $A_{\gamma,K}$ by
$$\sigma_T:=\sigma_{\psi_T,\delta_T}$$.
\end{ddd}

\begin{prop}\label{pos}
The functional $\sigma_T$ extends continuously to a
state $\sigma_T$ on $A_\gamma$.
\end{prop}
\proof
We choose a sequence $f_j\in C^K(K/M)$ such that
$|f_j|^2$ is a $\delta$-sequence located at $1M$.
In particular we require that $\|f_j\|_{L^2(K/M)}=1$ for all $j$.
Note that $\lim_j |f_j|^2 \psi_T=\delta_T$ and
$$\|f_j\psi_T\|_{L^2(K\times_MV_\kappa)\otimes V_\gamma}=1\ .$$

Let $f\in A_{\gamma,K}$.
We consider the finite subset $F=F(\psi_T,f)\subset \hat K$.
Then we can write 
\begin{eqnarray*}
\sigma_T(f)&=&
\sigma_{\psi_T,\delta_{T,F}}(f)\\&=&
\lim_j 
\sigma_{\psi_T,[|f_j|^2 \psi_T]_F} (f)\\&=&
\lim_j 
\sigma_{\psi_T,|f_j|^2 \psi_T} (f)\\&=&
\lim_j 
\sigma_{f_j\psi_T,f_j \psi_T} (f)
\end{eqnarray*}
using \ref{ds} in the last step. 
Since $\sigma_{f_j\psi_T,f_j \psi_T}$ is a  state
the assertion follows.

\subsubsection{}

Recall that $\sigma_T$ may depend on the choice of the sequence $\xi_n$.
Let us fix $\kappa\in \hat M$, $\gamma\in \hat K$, and
$T\in [V_\kappa\otimes V_\gamma]^M$ with $\|T\|=1$.
\begin{ddd}
For each regular $l\in L(\kappa)$ we define the
set $V(l,\kappa,\gamma,T)$ of states on $A_\gamma$ 
of the form $\sigma_T$ for the various $l$-conveniently arranged
sequences
$\xi_n$. These states are called representation theoretic lifts.
\end{ddd}

\subsubsection{}

Let $\sigma_T\in V(l,\kappa,\gamma,T)$ be associated to the
$l$-conveniently arranged sequence $\xi_n$.
We again consider the $\delta$-sequence $f_j\in C^K(K/M)$ as in the
proof of Proposition \ref{pos}.

\begin{prop}\label{ghk}
There exists a sequence of integers $n_j\to \infty$ such that
$\sigma^{\xi_{n_j}}_{f_{j}\psi_T,f_{j}\psi_T}$ weakly converges to
  $\sigma_T$
as $j\to\infty$.
\end{prop}
\proof
Let $f\in A_{\gamma,K}$. We have
$\sigma_T(f)=\lim_{n} \sigma^{\xi_n}_{\psi_T,\delta_{T,F}}(f)$,
where the finite subset $F:=F(\psi_T,f)\subset \hat K$ only depends on $T$ and $f$.
We estimate
\begin{eqnarray*}
| \sigma^{\xi_n}_{\psi_T,\delta_{T,F}}(f)- \sigma^{\xi_n}_{f_j\psi_T,f_j\psi_T}(f)|&\le&
| \sigma^{\xi_n}_{\psi_T,\delta_{T,F}}(f)-
 \sigma^{\xi_n}_{\psi_T,[|f_j|^2\psi_T]_F}(f)|\\&&+| \sigma^{\xi_n}_{\psi_T,|f_j|^2
\psi_T}(f)-
\sigma^{\xi_n}_{f_j\psi_T,f_j\psi_T}(f)|\ .
\end{eqnarray*}
We choose $n_j>j$ sufficiently large such that (by Corollary \ref{ds}) 
$$|\sigma^{\xi_{n_j}}_{\psi_T,|f_j|^2
\psi_T}(f)-
\sigma^{\xi_{n_j}}_{f_j\psi_T,f_j\psi_T}(f)|\le j^{-1}\ .$$
Since $\lim_j   [|f_j|^2\psi_T]_F = \delta_{T,F}$ inside a finite
dimensional vector space we
conclude that $$\lim| \sigma^{\xi_{n_j}}_{\psi_T,\delta_{T,F}}(f)-
\sigma^{\xi_{n_j}}_{f_j\psi_T,f_j\psi_T}(f)|=0\ .$$
 \hB

\subsubsection{}

If $\Gamma$ is arithmetic (this is automatic if $G$ has higher real
rank and $\Gamma\subset G$ is irreducible), then we can consider Hecke operators.
Let $H^\Gamma\subset H$ denote the subset of
$\Gamma$-invariant vectors in a representation $(\pi,H)$ of $G$.
If $h\in G$ is in the
commensurator of $\Gamma$, i.e. $\Gamma^h:=h\Gamma h^{-1}$ and $\Gamma$ are
commensurable, then we define the following operator
$T_h:H^\Gamma\rightarrow H^\Gamma$.
\begin{ddd}
$T_h(\phi)=\sum_{[\gamma]\in  \Gamma/(\Gamma\cap \Gamma^h)}
\pi(\gamma)\pi(h
)\phi$.
\end{ddd}

\subsubsection{}\label{fgh1}

We can apply the Hecke operator $T_h$  to $[H^*_{-\infty}]^\Gamma$
and to $C^\infty(\Gamma\backslash G)\otimes V_\gamma$.
Let $\eta\in \C$ and assume that  $\xi_n$ is a $l$-conveniently arranged
sequence of unitary embeddings $H^{\kappa,i\lambda_n}\rightarrow
L^2(\Gamma\backslash G)$ for some regular $l\in L(\kappa)$ such that $T_h\xi_n=\eta \xi_n$ for all $n$. 
Let $\sigma_T$ denote a limit state as in \ref{sit}.
For $u\in L^2(\Gamma\backslash G)\otimes V_\gamma$ let
$\sigma_u$ be the functional on $A_\gamma$ given by
$$\sigma_u(f)=\int_{\Gamma\backslash G} <u(\Gamma
g),f(\Gamma g)u(\Gamma g)>\ .$$

\begin{theorem}\label{dde}
There exists a sequence $u_j\in C^\infty(\Gamma\backslash G)\otimes V_\gamma$
of eigenvectors of $T_h$ to the eigenvalue $\eta$ such that
$\|u_j\|_{L^2(\Gamma\backslash G)\otimes V_\gamma}=1$ and
$\sigma_T$ is the weak limit of the sequence of states $\sigma_{u_j}$.
\end{theorem}
\proof
We choose the sequence $n_j$ as in Proposition \ref{ghk} and set
$$u_j:=\xi_{n_j}(f_j\psi_T)\ .$$
Since $\|f_j\psi_T\|_{L^2(K\times_MV_\kappa)\otimes V_\gamma}=1$ and
$\xi_{n_j}$ is an unitary embedding, the section $u_j$ is a unit
vector. It is a Hecke-eigenvector since $\xi_{n_j}$ is so.
\hB

\subsubsection{}\label{spaet}

Let $A_{\gamma}^K\subset A_\gamma$ be the subalgebra
of $K$-invariants, i.e. the algebra of sections of the bundle of
endomorphisms of the vector bundle $\Gamma\backslash G\times_K
V_\gamma\rightarrow \Gamma\backslash G/K$.
For each unit vector $u\in L^2(\Gamma\backslash G\times_K V_\gamma)$
we consider the state $\sigma_u$ as defined in \ref{fgh1}.

Let $\xi:H^{\kappa,i\lambda}\rightarrow L^2(\Gamma\backslash
G)$ be a
unitary embedding. 
Then we have a unit vector $\xi(\psi_T)\in
L^2(\Gamma\backslash G\times_K V_\gamma)=[L^2(\Gamma\backslash G)\otimes V_\gamma]^K$.
Let $f\in  A_{\gamma}^K$.
As in \ref{dif} the functional
$\phi\mapsto \sigma^{\xi}_{\psi_T,\phi}(f)$ extends to
distributions $\phi$.

\begin{lem}\label{mite}
We have 
$\sigma^{\xi}_{\psi_T,\delta_T}(f)=\sigma_{\xi(\psi_T)}(f)$.
\end{lem}
\proof
Let $F=F(\psi_T,f)\subset \hat K$
 be the finite subset of $K$-types as in \ref{dif}.
We may assume that $F$ contains the trivial $K$-type.
Note that
$$ \int_K\gamma(k)^{-1}\pi^{\kappa,i\lambda}(k^{-1})\delta_{T,F}=
 \int_K \gamma(k)^{-1}\pi^{\kappa,i\lambda}(k^{-1})\delta_{T}=
\psi_T \ .$$
In fact, the integral defines a $K$-invariant vector in
$H^{\kappa,i\lambda}\otimes V_\gamma$ which by Frobenius reciprocity
is determined by its value at $1\in K$. In view of the definition
of $\delta_T$ this evaluation is $T\in [V_\kappa\otimes V_\gamma]^M$.
We compute 
\begin{eqnarray*}
\sigma^{\xi}_{\psi_T,\delta_T}(f)&=&
\sigma^{\xi}_{\psi_T,\delta_{T,F}}(f)\\
&=&
\int_{\Gamma\backslash G} <\xi(\psi_T)(\Gamma
g),f(\Gamma g)\xi(\delta_{T,F})(\Gamma g)>\\
&=&\int_{\Gamma\backslash G} \int_K <\xi(\psi_T)(\Gamma
g),\gamma(k)f(\Gamma gk)\gamma(k)^{-1}\xi(\delta_{T,F})(\Gamma g)>\\
&=&\int_{\Gamma\backslash G} \int_K <\gamma(k^{-1})\xi(\psi_T)(\Gamma
gk^{-1}),f(\Gamma g)\gamma(k)^{-1}\xi(\delta_{T,F})(\Gamma gk^{-1})>\\
&=&\int_{\Gamma\backslash G} \int_K <\xi(\psi_T)(\Gamma
g),f(\Gamma
g)\gamma(k)^{-1}\xi(\pi^{\kappa,i\lambda}(k^{-1})\delta_{T,F})(\Gamma
g)>\\
&=&\int_{\Gamma\backslash G} \int_K <\xi(\psi_T)(\Gamma
g),f(\Gamma
g)\xi(\psi_T)(\Gamma
g)>\\
&=&\sigma^{\xi}_{\psi_T,\psi_T}(f)\\
&=&\sigma_{\xi(\psi_T)}(f)\ .
\end{eqnarray*}
\hB 

\subsubsection{}
Let $l\in L(\kappa)$ be regular and $\sigma_T\in V(l,\kappa,\gamma,T)$
associated to the $l$-conveniently arranged sequence of unitary embeddings
$\xi_n$. 
\begin{prop}\label{fr}
The sequence of states $\sigma_{\xi_n(\psi_T)}$ on 
$A_{\gamma}^K$ has a weak limit which is given by the restriction
of $\sigma_T$ from $A_\gamma$ to $A_{\gamma}^K$.
\end{prop}
\proof
This is a consequence of Lemma \ref{mite}.
\hB 

\section{Invariance of representation theoretic lifts}\label{dgl}

\subsubsection{}

Throughout this section we fix a unitary $K$-representation $(\gamma,V_\gamma)$
and an element $\kappa\in \hat M$. We further fix a regular $l\in
L(\kappa)$, and $T\in
[V_\kappa\otimes V_\gamma]^M$ such that $\|T\|=1$.

We want to
study the invariance properties of the states $\sigma\in
V(l,\kappa,\gamma,T)$ constructed in Section \ref{gf}.

\subsubsection{}

We first consider the action of the subgroup $M$. 
Note that $K$ acts on the $C^*$-algebra
$A_\gamma:=C(\Gamma\backslash G)\otimes \End(V_\gamma)$  by
\begin{equation}\label{act}
k(f\otimes A)(g):=(f\otimes \gamma(k)A\gamma(k^{-1}))(gk)\ ,\quad f\in C(\Gamma\backslash G), A\in \End(V_\gamma), k\in K\ .
\end{equation}
By duality, (\ref{act}) induces a $K$-action on the set of states of $A_\gamma$.

\subsubsection{} 

We consider $T$ as
an element of $\Hom_M(V_{\tilde \kappa},V_\gamma)$, where $\tilde
\kappa$ is the dual representation of $\kappa$. Then $\dim(V_\kappa)TT^*\in\End(V_\gamma)$ is an $M$-equivariant
projection. Let $P_T:=1\otimes \dim(V_\kappa)TT^*\in A_\gamma$ be the corresponding projection in $A_\gamma$.

\begin{prop}\label{mulm}
Let $\sigma\in V(l,\kappa,\gamma,T)$. 
\begin{enumerate}
\item $\sigma$ is $M$-invariant w.r.t. the action induced by (\ref{act}).
\item For each $f\in A_\gamma$ we have $\sigma(P_TfP_T)=\sigma(f)$.
\end{enumerate}
\end{prop}
\proof Let $\lambda\in\aaaa^*$, and let $\xi:H^{\kappa,i\lambda}\rightarrow L^2(\Gamma\backslash
G)$ be a unitary embedding. Using that
$$\psi_T,\delta_T\in [H_{-\infty}^{\kappa,i\lambda}\otimes V_\gamma]^M\ \mbox{ and }\  (\id\otimes \dim(V_\kappa)TT^*)\delta_T=\delta_T$$
we see that
$m\sigma^{\xi}_{\psi_T,\delta_T}=\sigma^{\xi}_{\psi_T,\delta_T}$ for all $m\in M$, and that $\sigma^{\xi}_{\psi_T,\delta_T}(P_Tf)=\sigma^{\xi}_{\psi_T,\delta_T}(f)$ for all
$f\in A_{\gamma,K}^\infty$. Taking the limit over an $l$-conveniently arranged sequence $\xi_n$
we obtain the first assertion and that
\begin{equation}\label{dick}
\sigma(P_Tf)=\sigma(f)\ .
\end{equation}
A state $\sigma$ on a $C^*$-algebra is a real functional, i.e., it satisfies
$\sigma(f^*)=\overline{\sigma(f)}$ for all $f$. Using (\ref{dick}) we obtain
$$ \sigma(P_TfP_T)=\sigma(fP_T)=\overline{\sigma(P_T f^*)}=\overline{\sigma(f^*)}=\sigma(f)\ .$$
Since $A_{\gamma,K}^\infty\subset A_\gamma$ is dense this finishes the proof of the proposition.
\hB

\subsubsection{}
The proposition tells us that $\sigma\in V(l,\kappa,\gamma,T)$ is given as a pull back of a state
$\sigma^0$ on the smaller algebra $P_TA_\gamma^MP_T$ which is in fact isomorphic to $C(\Gamma\backslash G\times_M\End(V_{\tilde\kappa}))$.

It seems to be likely
that $\sigma^0$ is actually a trace on $C(\Gamma\backslash G\times_M\End(V_{\tilde\kappa}))$. Up to now we do not know how to prove
this property. It would imply that $\sigma$ is determined by a
probability measure
on $\Gamma\backslash G/M$ alone. By Theorem \ref{inv} below this measure
would be right $A$-invariant.

\subsubsection{}

The right-regular representation of $G$ on $C(\Gamma\backslash G)$ induces
an action of $G$ by automorphisms on the $C^*$-algebra
$A_\gamma=C(\Gamma\backslash G)\otimes \End(V_\gamma)$. Dually, we obtain
a $G$-action on the states of $A_\gamma$. The main goal of this section
is to prove the following higher rank analog of Proposition
\ref{inv1}.
Recall the Iwasawa decomposition $G=KAN$ (see \ref{hs}).

\begin{theorem}\label{inv}
The states $\sigma\in V(l,\kappa,\gamma,T)$ are $A$-invariant.
\end{theorem}

Let $\xi_n:H^{\kappa,i\lambda_n}\rightarrow L^2(\Gamma\backslash
G)$ be a $l$-conveniently arranged sequence of unitary embeddings giving rise to $\sigma$.
Following the approach of \cite{sv04}, Section 4, we will
exhibit a certain family of differential operators $D(\lambda)$,
depending polynomially on
$\lambda\in\aaaa^*$, such that
$D(\lambda_n)\sigma^{\xi_n}_{\psi_T,\delta_T}=0$. In the limit
$n\to\infty$ this will imply $A$-invariance of $\sigma$.

\subsubsection{}\label{cent}

For any real or complex Lie algebra $\laaa$ let $\cU(\laaa)$ be its
universal enveloping algebra over $\C$. The algebra $\cD_\gamma:=\cug\otimes\End(V_\gamma)\otimes\End(V_\gamma)^\opp$ acts by differential
operators on $A_\gamma^\infty=C^\infty(\Gamma\backslash G)\otimes \End(V_\gamma)$:
$$(X\otimes A\otimes B)(f\otimes C):=Xf\otimes ABC\ ,\quad X\in\cug,\ A,B,C\in\End(V_\gamma),\ f\in C^\infty(\Gamma\backslash G)\ .$$
The subspace of $K$-finite elements
$A_{\gamma,K}^\infty:=A_\gamma^\infty\cap A_{\gamma,K}$ is invariant w.r.t. this action.
Therefore we have an action of $\cD_\gamma$ on the space
of functionals on $A_{\gamma,K}^\infty$ given by
\begin{equation}\label{kuddel}
D\sigma(f):=\sigma(D^{\, t} f)\ ,
\end{equation}
where $D\mapsto D^{\, t}$ is the anti-automorphism of $\cD_\gamma$
induced by
$$(X\otimes 1)^t=-X\otimes 1\otimes 1,\  X\in\gaaa,\quad (1\otimes A\otimes B)^t=1\otimes B\otimes A,\ A,B\in\End(V_\gamma)\ . $$

We are mainly concerned with the subalgebra
$$\cL_\gamma:=\cU(\naaa\oplus\aaaa)\otimes\End(V_\gamma)^\opp
\subset\cug\otimes\End(V_\gamma)\otimes\End(V_\gamma)^\opp=\cD_\gamma\ .$$
There is a linear map $q_\gamma: \cug\rightarrow\cL_\gamma$ which sends
$X\otimes Y\in \cU(\naaa\oplus\aaaa)\otimes \cU(\kaaa)\cong\cug$ to
${X\otimes\gamma(Y)}\in \cU(\naaa\oplus\aaaa)\otimes\End(V_\gamma)^\opp=\cL_\gamma$.

If $\laaa$ is a Lie algebra and $\vp:\laaa\rightarrow \C$ is a Lie
algebra
homomorphism, then there is a corresponding translation automorphism
$\tau_\vp:\cU(\laaa)\rightarrow\cU(\laaa)$ characterized by $\tau_\vp(X)=X+\vp(X)\cdot 1$, $X\in\laaa$.

\subsubsection{}

Let $\maaa$ be the Lie algebra of $M$. We choose a Cartan subalgebra $\taaa\subset\maaa$. Then $\haaa:=\taaa\oplus\aaaa$ is a Cartan subalgebra
of $\gaaa$. Let $W_\C$ be the Weyl group of $\haaa_\C$ in $\gaaa_\C$.
Let $P\in\cU(\haaa)^{W_\C}$. We view $P$ as a complex-valued polynomial
on $\haaa_\C^*=\taaa^*_\C\oplus\aaaa_\C^*$.
Its differential $P'$ is a polynomial
on $\haaa_\C^*$ with values in $\haaa_\C\cong(\haaa_\C^*)^*$.
Let $\mu_\kappa\in i\taaa^*\subset\haaa_\C^*$ be an extremal weight
of $\kappa$, i.e., it is the highest weight of $\kappa$ w.r.t. some
positive root system of $\taaa$ in $\maaa_\C$.

\begin{prop}\label{main}
Fix $P\in\cU(\haaa)^{W_\C}$ of degree $\le d\in\nat_0$. Then there exists
a $\cL_\gamma$-valued polynomial $J_P$ on $\aaaa_\C^*$ of degree at most $d-2$ such that
for all unitary $G$-maps $\xi:H^{\kappa,i\lambda}\rightarrow L^2(\Gamma\backslash
G)$
$$ \left(q_\gamma\circ\tau_{\mu_\kappa}(P'(i\lambda))+J_P(i\lambda)\right)\sigma^{\xi}_{\psi_T,\delta_T}=0\ .$$
Here $P'(i\lambda)\in\haaa_\C$ is viewed as an element of $\cU(\haaa)$.
Then $\tau_{\mu_\kappa}(P'(i\lambda))\in \cU(\haaa)\subset\cug$, and $q_\gamma$ can be applied.
\end{prop}

The Weyl group $W_0$ of $\taaa_\C$ in $\maaa_\C$ considered as subgroup of $W_\C$ fixes the element $i\lambda\in\haaa_\C^*$. It follows that $P'(i\lambda)\in\haaa_\C^{W_0}$. Since all extremal weights of $\kappa$ are conjugated
by $W_0$, the element $\tau_{\mu_\kappa}(P'(i\lambda))\in \cU(\haaa)$
does not depend on the choice of $\mu_\kappa$.

\subsubsection{}

The proof of Proposition \ref{main} starts in the next paragraph
\ref{next2} and will then occupy the remainder of this section.
Here we argue as in \cite{sv04}, Corollary 4.6, Lemma 4.7 and Corollary 4.8
in order to conclude that Proposition \ref{main} implies Theorem \ref{inv}.

We first assume $P\in\cU(\haaa)^{W_\C}$ to be homogeneous of degree
$d$.  Fix $f\in
A_{\gamma,K}^\infty$. Then 
 by Proposition \ref{main} the equation
$$ \sigma^{\xi_n}_{\psi_T,\delta_T}\left(\left(q_\gamma\circ\tau_{\mu_\kappa}
\left(P'\left(\frac{i\lambda_n}{\|\lambda_n\|}\right)\right)
+\frac{J_P(i\lambda_n)}{\|\lambda_n\|^{d-1}}\right)^{t}f\right)=0\ $$
holds for all $n$. There is a finite dimensional subspace $V\subset A_{\gamma,K}^\infty$ such that
$$f_n:=\left(q_\gamma\circ\tau_{\mu_\kappa}
\left(P'\left(\frac{i\lambda_n}{\|\lambda_n\|}\right)\right)
+\frac{J_P(i\lambda_n)}{\|\lambda_n\|^{d-1}}\right)^{t}f\in V$$
for all $n$. Moreover, $f_n$ converges in $V$
to $\left(q_\gamma\circ\tau_{\mu_\kappa}(P'(il))\right)^{\,t}(f)$. We obtain
$$\left(q_\gamma\circ\tau_{\mu_\kappa}(P'(il))\sigma_T\right)(f)=0\ .$$
Since this is valid for all $f$, and each $P\in\cU(\haaa)^{W_\C}$ can be decomposed into homogeneous components, we conclude that $\sigma$ is annihilated
by all the operators $q_\gamma\circ\tau_{\mu_\kappa}(P'(il))$, $P\in\cU(\haaa)^{W_\C}$.

For each $H\in\aaaa$ there exists an element $P_H\in\cU(\haaa)^{W_\C}$
such that $P_H'(il)=H$ (see \cite{sv04}, Lemma 4.7). Here the regularity of
$l$ is crucial. Then $q_\gamma\circ\tau_{\mu_\kappa}(P_H'(il))=H\otimes 1$.
It follows that $\sigma$ is $\aaaa$-invariant, and hence $A$-invariant.
This proves Theorem \ref{inv} assuming Proposition \ref{main}. \hB 

\subsubsection{}\label{next2}

Let $\tilde\kappa$, $\tilde\gamma$ be the representations dual to $\kappa$, $\gamma$.
Let $\nu\in\aaaa_\C^*$. The algebra $\cug$ acts on the tensor product
representation $H^{\tilde\kappa,-\nu}_{\infty}\otimes H^{\kappa,\nu}_{-\infty}$.
We let $A\otimes B\in\End(V_{\gamma})\otimes\End(V_{\gamma})^\opp$ act on $V_{\tilde\gamma}\otimes V_\gamma$
by $B^\top\otimes A$, where
$B^\top\in \End(V_{\tilde\gamma})$ is the dual operator
of $B\in\End(V_\gamma)$.
We obtain an action of $\cD_\gamma$ on
$(H^{\tilde\kappa,-\nu}_\infty\otimes V_{\tilde\gamma}) \otimes (H^{\kappa,\nu}_{-\infty}\otimes V_\gamma)$. If $\nu\in i\aaaa^*$, then we have a canonical antilinear identification $R: H^{\kappa,\nu}\otimes V_\gamma\rightarrow
H^{\tilde\kappa,-\nu}\otimes V_{\tilde\gamma}$.
It is a direct consequence of Lemma \ref{leib} and (\ref{kuddel})
that,
\begin{equation}\label{muddel}
\mbox{if }\ D(R\psi\otimes \phi)=\sum_i \psi_i\otimes \phi_i\, ,\quad
\mbox{ then }\ D\sigma^\xi_{\psi,\phi}=\sum_i\sigma^\xi_{R^{-1}\psi_i,\phi_i}\, .
\end{equation}
Here $\xi:H^{\kappa,\nu}\rightarrow L^2(\Gamma\backslash
G)$, $\psi\in H^{\kappa,\nu}_\infty\otimes V_\gamma$, $\phi\in H^{\kappa,\nu}_{-\infty}\otimes V_\gamma$, $D\in\cD_\gamma$.

\subsubsection{}

The composition of  $\nu\in\aaaa_\C^*$ with the projection of
$\naaa\oplus\aaaa\rightarrow \aaaa$ defines a Lie algebra homomorphism $\nu:\naaa\oplus\aaaa\rightarrow \C$.
Let $\tau_\nu$ be the corresponding translation automorphism of $\cU(\naaa\oplus\aaaa)$ (see \ref{cent}). Then $\tau_\nu\otimes \id_{\End(V_{\gamma})^\opp}$
is an automorphism of $\cL_\gamma$ which will be denoted by $\tau_\nu$ as well.

\begin{lem}\label{engels}
For $\nu\in \aaaa_\C^*$, $\psi\in H^{\tilde\kappa,-\nu}_\infty$,
$T\in [V_\kappa\otimes V_\gamma]^M$, and $D\in\cL_\gamma$
we have
$$\tau_{\nu+\rho}(D)(\psi\otimes\delta_T)=(D\psi)\otimes\delta_T\ .$$
\end{lem}
\proof
It suffices to check the assertion for the generators $X\in\aaaa$, $Y\in\naaa$,
and $B\in\End(V_\gamma)^\opp$ of $\cL_\gamma$. For $D=B$ the assertion holds
by definition while for $D=Y$ it follows from $Y\delta_T=0$.
Now let $X\in\aaaa$ and $\phi\in H^{\tilde\kappa,-\nu}_\infty$.  Then
$$\langle X\delta_T,\phi \rangle =-\langle \delta_T,X\phi\rangle =-\langle T,X\phi(1)\rangle=-\langle T,(\nu+\rho)(X)\phi(1)\rangle=-(\nu+\rho)(X)\langle \delta_T,\phi \rangle\ .$$
Hence $X\delta_T=-(\nu+\rho)(X)\delta_T$ and
$$\tau_{\nu+\rho}(X)(\psi\otimes\delta_T)=(X\psi)\otimes\delta_T +
\psi\otimes (X\delta_T)+(\nu+\rho)(X)(\psi\otimes\delta_T)=(X\psi)\otimes\delta_T\ .$$
\hB

\subsubsection{}

Let $\uaaa$ ($\bar\uaaa$, resp.) be the sum of positive (negative) root
spaces in $\gaaa_\C$ w.r.t. a chosen positive Weyl chamber in $\haaa_\R:=i\taaa\oplus \aaaa$. We arrange this choice such 
that $\naaa_\C\subset \uaaa$.
Then we have a decomposition
$$ \gaaa_\C=\uaaa\oplus \haaa_\C\oplus \bar\uaaa\ . $$
It induces decomposition
$$ \cug =\cU(\haaa)\oplus (\uaaa\cug +\cug\bar\uaaa)\ .$$
Let $p$ be the projection onto the first summand. Note that $p$ is equivariant
w.r.t. the adjoint action of $\haaa_\C$ on $\cug$. By $\cU^{\le d}(\laaa)$ we
denote the subspace of elements of $\cU(\laaa)$ of degree at most $d$.
The following lemma is essentially Lemma 4.3 of \cite{sv04}.

\begin{lem}\label{murks}
If $Z\in \cU^{\le d}(\gaaa)^{\haaa_\C}$, then
$Z-p(Z)\in\cU(\naaa)\cU^{\le d-2}(\aaaa)\cU(\kaaa)$.
\end{lem}
\proof Observe that
$$(\uaaa\cug +\cug\bar\uaaa)^{\haaa_\C}\subset \uaaa\cug\bar\uaaa\ .$$
If $Z\in \cU^{\le d}(\gaaa)^{\haaa_\C}$, then by $\haaa_\C$-equivariance of $p$
\begin{equation}\label{arm}
Z-p(Z)\in\uaaa\cU^{\le d-2}(\gaaa)\bar\uaaa\subset \cU(\uaaa)\cU^{\le d-2}(\haaa)
\cU(\bar\uaaa)\ .
\end{equation}
Using that $\haaa_\C\subset\aaaa_\C\oplus\maaa_\C$, $\uaaa\subset \naaa_\C\oplus \maaa_\C$, $\maaa_\C\subset\kaaa_\C$, and $\bar\uaaa\subset \naaa_\C\oplus\kaaa_\C$ one shows inductively that
the right hand side of (\ref{arm}) is contained in $\cU(\naaa)\cU^{\le
  d-2}(\aaaa)\cU(\kaaa)$. This proves the lemma.
\hB

\subsubsection{}

Let $\cZ(\gaaa)$ be the center of $\cug$. 
If $\nu\in\aaaa_\C^*$, then $\cZ(\gaaa)$ acts on $H^{\tilde\kappa,-\nu}_\infty$
by a certain character denoted by
$\chi_{\kappa,\nu}$. Recall the definition of $q_\gamma:\cug\rightarrow\cL_\gamma$ from \ref{cent}.
For $Z\in\cZ(\gaaa)$ we consider the elements $p_\gamma(Z):=q_\gamma(p(Z))\in\cL_\gamma$ and  $b_\gamma(Z):=q_\gamma(Z-p(Z))\in\cL_\gamma$.

\begin{lem}\label{marx}
If $\psi\in [H^{\tilde\kappa,-\nu}_\infty\otimes V_{\tilde\gamma}]^K$, then we have for all $Z\in\cZ(\gaaa)$
$$ \left(p_\gamma(Z)-\chi_{\kappa,\nu}(Z)+b_\gamma(Z)\right)\psi=0\ .$$
\end{lem}
\proof
Let $W\in\cug$, $Y\in\kaaa$ and $\psi\in
[H^{\tilde\kappa,-\nu}_\infty\otimes V_{\tilde\gamma}]^K$. Then we have
$q_\gamma(WY)=q_\gamma(Y)q_\gamma(W)$ and
$$(WY\otimes 1)\psi=(W\otimes \tilde\gamma(-Y))\psi=(W\otimes \gamma(Y)^\top)\psi=\left(q_\gamma(Y)(W\otimes 1)\right)\psi\ .$$
It follows by induction that
$ (X\otimes 1)\psi=q_\gamma(X)\psi$
for any $X\in\cug\cong \cU(\naaa\oplus\aaaa)\cU(\kaaa)$.
Applying this to $X=Z\in\cZ(\gaaa)$ we obtain
$$0= \left((Z\otimes 1)-\chi_{\kappa,\nu}(Z)\right)\psi=\left(q_\gamma(Z)-\chi_{\kappa,\nu}(Z)\right)
\psi=\left(p_\gamma(Z)+b_\gamma(Z)-\chi_{\kappa,\nu}(Z)\right)\psi\ .$$
\hB
Combining Lemma \ref{marx} with Lemma \ref{engels} we obtain

\begin{kor}\label{fast}
If $\psi\in [H^{\tilde\kappa,-\nu}_\infty\otimes V_{\tilde\gamma}]^K$ and
$T\in [V_\kappa\otimes V_\gamma]^M$, then we have for all $Z\in\cZ(\gaaa)$
$$ \left(\tau_{\nu+\rho}(p_\gamma(Z))-\chi_{\kappa,\nu}(Z)
+\tau_{\nu+\rho}(b_\gamma(Z))\right)
\left(\psi\otimes\delta_T\right)=0\ .$$
\end{kor}

\subsubsection{}

Let $\rho_\haaa\in\haaa_\C^*$ be given by $\rho_\haaa(H)=\frac{1}{2}\Tr\, \Ad(H)_{|\uaaa}$. Then $\rho_\taaa:=\rho_\haaa-\rho\in i\taaa^*\subset \haaa_\C^*$. The composition  $\tau_{\rho_\haaa}\circ p$ maps
the algebra
$\cZ(\gaaa)$ isomorphically onto $\cU(\haaa)^{W_\C}$. In fact, this
map is the celebrated
Harish-Chandra isomorphism. If $P\in\cU(\haaa)^{W_\C}$, then we denote its
preimage under the Harish-Chandra isomorphism by $Z_P$. The roots of
$\taaa_\C$ in $\maaa_\C\cap\uaaa$ form a positive root system of $\taaa_\C$
in $\maaa_\C$. Let $\mu_\kappa\in i\taaa^*$ be the highest weight of $\kappa$
with respect to this system of positive roots.

\begin{lem}\label{trans}
Fix $P\in\cU(\haaa)^{W_\C}$ of degree $\le d$. Then the $\cL_\gamma$-valued polynomial
$J_P$ on $\aaaa_\C^*$ defined by
$$ J_P(\nu):=\tau_{\nu+\rho}(p_\gamma(Z_P))-\chi_{\kappa,\nu}(Z_P)
+\tau_{\nu+\rho}(b_\gamma(Z_P))-q_\gamma\circ\tau_{\mu_\kappa}(P'(\nu))$$
has degree at most $d-2$. Here the expression $q_\gamma\circ\tau_{\mu_\kappa}(P'(\nu))$ is interpreted as in Proposition
\ref{main}.
\end{lem}
\proof If $P$ has degree $d$, then $Z_P\in \cU(\gaaa)^{\le d}$. Now Lemma \ref{murks} implies that $\nu\mapsto\tau_{\nu+\rho}(b_\gamma(Z_P))$ has degree at most $d-2$.

We now analyze the first two terms appearing in the definition of $J_P$.
It is well-known that $\chi_{\kappa,\nu}(Z_P)=P(\nu-\rho_\taaa-\mu_\kappa)$.
The reader may verify this by computing the value $\langle w, Z_P\psi(1)\rangle$,
where $w\in V_\kappa$ is the highest weight vector and $\psi\in H^{\tilde\kappa,-\nu}$ (compare the proof of Lemma \ref{engels}). Let $S\in\cU(\haaa)$. Choosing a basis of $V_\gamma$ consisting of weight vectors w.r.t. the action of $\taaa$ we may view $q_\gamma(S)\in\cU(\aaaa)\otimes\End(V_\gamma)^\opp$
as a diagonal matrix whose entries are polynomials on $\aaaa_\C^*$.
The matrix entry corresponding to a weight $\mu_i\in i\taaa^*$ is then given
by the polynomial ${\aaaa_\C^*\ni x}\mapsto S(x+\mu_i)$.
Therefore the matrix entries of $p_\gamma(Z_P)=q_\gamma(p(Z_P))$ are
$\aaaa_\C^*\ni x\mapsto P(x-\rho_\haaa+\mu_i)$, and the ones of
$\tau_{\nu+\rho}(p_\gamma(Z_P))-\chi_{\kappa,\nu}(Z_P)$ are given by $x\mapsto
P(x+\nu-\rho_\taaa+\mu_i)-P(\nu-\rho_\taaa-\mu_\kappa)$. The
Taylor expansion of $P$ at $\nu-\rho_\taaa-\mu_\kappa$ yields
$$P(x+\nu-\rho_\taaa+\mu_i)-P(\nu-\rho_\taaa-\mu_\kappa)
=P'(\nu-\rho_\taaa-\mu_\kappa)(x+\mu_i+\mu_\kappa)
+Q_\nu(x+\mu_i+\mu_\kappa)\ ,$$
where $P'(\nu-\rho_\taaa-\mu_\kappa)$ is viewed as a linear form on $\haaa_\C^*$ and the polynomial $Q_\nu$ is formed by partial derivatives
of $P$ at $\nu-\rho_\taaa-\mu_\kappa$ of degree at least $2$. Therefore the degree
of $Q_\nu$ w.r.t. $\nu$ is bounded by $d-2$.
Using that the degree of $\nu \mapsto P'(\nu-\rho_\taaa-\mu_\kappa)-P'(\nu)$
is also bounded by $d-2$ we conclude that the same is true for
$$\nu\mapsto P(x+\nu-\rho_\taaa+\mu_i)-P(\nu-\rho_\taaa-\mu_\kappa)
-P'(\nu)(x+\mu_i+\mu_\kappa)\ .$$
It follows that the degree of the $\cL_\gamma$-valued polynomial
$$\nu\mapsto\tau_{\nu+\rho}(p_\gamma(Z_P))-\chi_{\kappa,\nu}(Z_P)-q_\gamma\circ\tau_{\mu_\kappa}(P'(\nu))$$
is at most $d-2$. Since we have already estimated the degree of $\nu\mapsto\tau_{\nu+\rho}(b_\gamma(Z_P))$ the proof of the lemma is now complete.
\hB

\subsubsection{}

It is now easy to finish the proof of Proposition \ref{main}. Let $J_P$ be as
in Lemma \ref{trans}. Put $\lambda\in\aaaa^*$ and $T\in [V_\kappa\otimes V_\gamma]^K$. We form the corresponding elements $\psi_T$, $\delta_T$ in $H^{\kappa,i\lambda}_{\pm\infty}\otimes V_\gamma$. Then we have by Corollary \ref{fast}
$$ \left(q_\gamma\circ\tau_{\mu_\kappa}(P'(i\lambda))+J_P(i\lambda)\right)
(R(\psi_T)\otimes\delta_T)=0\ .$$
Now we apply formula (\ref{muddel}).

\section{The relation of the microlocal and representation theoretic
  lifts}\label{luzio}

\subsubsection{}

The discussion in the present section is completely parallel to
\cite{sv04}, Sec. 5.4. Its goal is to provide the link between the microlocal
and the representation theoretic lifts. The main result is
Corollary \ref{rg3}.

\subsubsection{}

Let $\gaaa=\kaaa\oplus\paaa$ be the Cartan decomposition.
Then we can write the cotangent bundle of $M:=\Gamma\backslash G/K$
as $T^*M=\Gamma\backslash G\times_K \paaa^*$.
We let $\pi:\Gamma\backslash G\times \paaa^*\rightarrow T^*M$ denote the
projection. 

The Killing form of $\gaaa$ restricts to a metric on $\paaa$ which is
$K$-invariant. It induces a Riemannian metrics on $TM$ and $T^*M$.
We further get an orthogonal decomposition $\paaa=\aaaa\oplus \aaaa^\perp$.
This induces an embedding $\aaaa^*\hookrightarrow \paaa^*$.

Let $SM\subset T^*M$ denote the unit cosphere bundle. Then $\pi$
restricts to a map $$q:\Gamma\backslash G\times S(\aaaa^*)\rightarrow SM\ .$$

\subsubsection{}

We now consider a unitary representation $\gamma\in \hat K$. It gives
rise to a bundle
$V(\gamma):=\Gamma\backslash G\times_KV_\gamma$ over $M$.
Let $p:SM\rightarrow M$ be the projection. We consider the identification
$$\Gamma\backslash G\times S(\aaaa^*)\times V_\gamma\rightarrow q^*\circ p^* V(\gamma)$$
defined such that $(\Gamma g,\lambda,v)$ corresponds to the point
$[\Gamma g,v]\in V(\gamma)$ in the fibre of $q^*\circ
p^* V(\gamma) $ over $(\Gamma g,\lambda)$.

In a similar manner we obtain an identification
$$\Gamma\backslash G\times S(\aaaa^*)\times \End(V_\gamma)\rightarrow
q^*\circ p^* \End(V(\gamma))\ .$$

\subsubsection{}

Let $f\in C(SM,p^*\End(V_\gamma))$. Then
$q^*f\in C(\Gamma\backslash G\times S(\aaaa^*))\otimes
\End(V_\gamma)$.
For $\lambda\in S(\aaaa^*)$ we define
$\tilde f_\lambda\in A_\gamma$ to be the restriction of $q^*f$
to $\Gamma\backslash G\times \{\lambda\}$.
 The map
$f\mapsto \tilde f_\lambda$ is a homomorphism of $C^*$-algebras
$$I_\lambda:C(SM,p^*\End(V_\gamma))\rightarrow A_\gamma\ .$$

\subsubsection{}

We now consider 
a regular element $l\in L(\kappa)\subset S(\aaaa^*)$ and $T\in [V_\kappa\otimes V_\gamma]^M$ with
$\|T\|=1$. Let $\xi_n$ be an $l$-conveniently arranged sequence giving
rise to a representation theoretic lift
$\sigma\in V(l,\kappa,\gamma,T)$.

We have a sequence of normalized vectors
$\xi_n(\psi_T)\in L^2(\Gamma\backslash G\times_K V_\gamma)$.
These sections are in fact smooth.
They give rise to functionals
$\sigma_{\xi_n(\psi_T)}$ on $C(SM,p^*\End(V_\gamma))$
(see \ref{statdef}).
After taking a subsequence we can and will assume that the sequence
$\sigma_{\xi_n(\psi_T)}$ converges weakly to some limit state, which we denote
by $\sigma_{micro}$ here. It is the microlocal lift
associated with the family of eigensections $\xi_n(\psi_T)\in L^2(\Gamma\backslash G\times_K V_\gamma)$ considered in Section \ref{gen}.

On the other hand we have functionals 
$\sigma^{\xi_n}_{\psi_T,\delta_T}$ on $A_{\gamma,K}$
defined in (\ref{sswe}). In fact, the same discussion as in \ref{dif}
shows that for $K$-finite $f$ and $\phi$ one can extend
$\sigma^{\xi}_{\phi,\psi}(f)$ to distributions $\psi$.

\subsubsection{}

Let $o(1)$ denote a quantity which tends to zero as $n$ tends to
infinity.

\begin{theorem}\label{hgf}
Assume that $f\in C(SM,p^*\End(V_\gamma))$ is such that $\tilde
f_\lambda$ is $K$-finite. Then we have  
$$\sigma_{\xi_n(\psi_T)}(f)=\sigma^{\xi_n}_{\psi_T,\delta_T}(\tilde
f_l)+o(1)\ .$$
\end{theorem}

\begin{kor}\label{rg3}
We have $\sigma_{micro}=I_l^*(\sigma)$.
In particular, $\sigma_{micro}$ is supported on $q(\Gamma\backslash G\times\{l\})$.
\end{kor}

\subsubsection{}

The idea of the proof is to verify the theorem on the symbols of
pseudodifferential operators $D(d,U,b)$ defined below.  This is the
contents of Proposition \ref{pp8}.
Then we show $\sigma_{micro}$ is supported on  $q(\Gamma\backslash
G\times\{l\})$ (Lemma \ref{ssk}). 
Finally we use that
these symbols span a dense subspace of $C(q(\Gamma\backslash
G\times\{l\}),p^*\End(V_\gamma))$ (Lemma \ref{denss}).

\subsubsection{}
 
We start with the construction of the family $D(d,U,b)$ of zero order
pseudodifferential operators  on $C^\infty(M,V(\gamma))$,
where $U\in \cU(\gaaa)^{\le d}$ and $b\in C^{\infty.K}(\Gamma\backslash
G)\otimes \End(V_\gamma)$.
Let $\Omega\in \cZ(\gaaa)$ be the Casimir operator. Note that
$\Omega+i$ is invertible on $L^2(M,V(\gamma))$ so that we can define
$(\Omega+i)^{-d/2}$ by the spectral theorem.
We identify
$L^2(M,V(\gamma))$ with the subspace of $K$-invariants
$[L^2(\Gamma\backslash G)\otimes V_\gamma]^K$.
Let $I_K:L^2(\Gamma\backslash G)\otimes V_\gamma\rightarrow
L^2(M,V(\gamma))$ denote the orthogonal projection.
It is given by
$$I_K(f)(\Gamma g)=\int_{K} \gamma(k) f(\Gamma gk)\ .$$
Then we define
$$D(d,U,b):=I_K\circ M_b\circ (R(U)\otimes 1)\circ (\Omega+i)^{-d/2}\ .$$
Here $R(U)$ denotes the right-regular action of $U$ on
$C^\infty(\Gamma\backslash G)$, and $M_b$ is the multiplication by
$b$, i.e.
$(M_bf)\Gamma g)=b(\Gamma g) f(\Gamma g)$.

The composition $I_K\circ M_b\circ (R(U)\otimes 1)$ is a differential
operator of order $\le d$. Since $(\Omega+i)^{-d/2}$
is a pseudodifferential operator of order $-d$ we conclude that
$D(d,U,b)$ is a pseudodifferential operator of order zero.

\subsubsection{}\label{lj1}

In this paragraph we compute the symbol of $D(d,U,b)$.
We consider the symmetrization map $sym:Sym(\paaa)\rightarrow U(\gaaa)$
defined on the degree $r$-subspace  by
$$sym(X_1\otimes \dots\otimes X_r)=\frac{1}{r!}\sum_{\sigma\in S_r}
X_{\sigma(1)}\dots X_{\sigma(r)}\ .$$
It extends to an isomorphism of vector spaces
$$\Phi:=mult(sym\otimes \id):Sym(\paaa)\otimes U(\kaaa)\rightarrow
U(\gaaa)\ .$$
This map preserves the filtrations by degree on both sides.
There exists a uniquely determined
$u\in S(\paaa)^d$ and $r\in S(\paaa)^{\le d-1}\otimes U(\kaaa)$ such
that $\Phi(u\otimes 1+r)=U$. Note that $U-\Phi(u\otimes 1)$
acts as a differential operator of order $\le d-1$.

We can now compute the symbol of $I_K\circ M_b\circ (R(U)\otimes
1)$. At the point
$[\Gamma g,\lambda]\in SM=\Gamma\backslash G\times_K S(\paaa^*)$
it is given by
$$s(D(d,U,b)(\Gamma g,\lambda))=\int_K u(\lambda^{k^{-1}})
\gamma(k)b(\Gamma g k)\gamma(k)^{-1} \in \End(V_\gamma)\ . $$
Note that $\widetilde{s(D(d,U,b))}_\lambda$ is $K$-finite. 
 
\begin{prop}\label{pp8}
The assertion of Theorem \ref{hgf} holds true for the functions of the
form
$s(D(d,U,b))\in C^\infty(SM,p^*\End(V_\gamma))$.
\end{prop}

\subsubsection{}

We know already that
$$<\xi_n(\psi_T),D(d,U,b)\xi_n(\psi_T)>=\sigma_{\xi_n(\psi_T)}\big(s(D(d,U,b))\big)+o(1)\
.$$
In order to prove the proposition we rewrite the left hand side.
Note that $$\xi_n(\psi_T)(\Gamma g\Omega)=(\lambda_n^2
+c(\kappa))\xi_n(\psi_T)(\Gamma g)\ ,$$
where $\xi_n\in H^{\kappa,i\lambda_n}_{-\infty}$ determines
$\lambda_n$, and $c(\kappa)\in \C$
is some constant independent of $n$. Therefore
we have
$$(\Omega+i)^{-d/2}\xi_n(\psi_T)=(\lambda_n^2+c(\kappa)+i)^{-d/2}\xi_n(\psi_T)\
.$$
%We get
%$$
%<\xi_n(\psi_T),D(d,U,b)\xi_n(\psi_T)>=
%(\lambda_n^2+c(\kappa)+i)^{-d/2}\int_{\Gamma\backslash G} 
%<\xi_n(\psi_T)(\Gamma g), b(\Gamma g)\xi_n(\psi_T)(\Gamma g U)>
%\ .$$
We thus get
\begin{eqnarray*}\lefteqn{
<\xi_n(\psi_T),D(d,U,b)\xi_n(\psi_T)>}&&\\
&=&
(\lambda_n^2+c(\kappa)+i)^{-d/2}\int_{\Gamma\backslash G}\int_K 
<\xi_n(\psi_T)(\Gamma g), \gamma(k)b(\Gamma gk)\xi_n(\psi_T)(\Gamma gk
U)>\\
&=&
(\lambda_n^2+c(\kappa)+i)^{-d/2}\int_{\Gamma\backslash G}\int_K
<\xi_n(\psi_T)(\Gamma g), \gamma(k) b(\Gamma g k)\gamma(k)^{-1}
\xi_n\left(\pi(U^k)\otimes 1)\psi_T\right)(\Gamma g)>\ ,
\end{eqnarray*}
where $U^k:=\Ad(k)(U)$ and $\pi:=\pi^{\kappa,i\lambda_n}$.

\subsubsection{}

We now use that
$\int_K(\pi(k)\otimes \gamma(k))\delta_T=\psi_T$ in order to write
\begin{eqnarray*}
(\pi(U^k)\otimes 1)
  \psi_T&=&(\pi(U^k)\otimes 1) \int_K (\pi(h)\otimes
  \gamma(h))\delta_T\\
&=&\int_K (\pi(h)\otimes \gamma(h)) (\pi(U^{h^{-1}k})\otimes
1)\delta_T
\end{eqnarray*}
Now we use that 
$$(\pi(U^k)\otimes 1)\delta_T=
u(\lambda^{k^{-1}}_n) \delta_T+\|\lambda_n\|^do(1)\ .$$
It follows that 
$$(\pi(U^k)\otimes 1)\psi_T=\int_K u(\lambda_n^{k^{-1}h}) 
(\pi(h)\otimes \gamma(h))\delta_T +\|\lambda_n\|^d o(1)\ .$$  
Our final rewriting is
\begin{eqnarray}\label{eq23}\lefteqn{
<\xi_n(\psi_T),D(d,U,b)\xi_n(\psi_T)>}&&\\&=&\int_{\Gamma\backslash G}\int_K\int_K
  <\xi_n(\psi_T)(\Gamma g), \gamma(k)^{-1}b(\Gamma g
k^{-1})\gamma(k)  u(l^{k^{-1}h})  \xi_n((\pi(h)\otimes \gamma(h))\delta_T)(\Gamma g)>
+o(1) \nonumber\ .\end{eqnarray}

\subsubsection{}

We now consider the right-hand side  of the equation \ref{hgf}.
We have
\begin{eqnarray}
\lefteqn{\sigma^{\xi_n}_{\psi_T,\delta_T}(\widetilde
D(d,U,b)_l)}&&\nonumber\\
&=&\int_{\Gamma\backslash G} <\xi_n(\psi_T)(\Gamma g), s(D,U,b)(\Gamma
g,l) \xi_n(\delta_T)(\Gamma g)>\nonumber\\
&=&\int_{\Gamma\backslash G} \int_K \int_K <\gamma(k)\xi_n(\psi_T)(\Gamma
gk),  \gamma(h)b(\Gamma gh)\gamma(h)^{-1}u(l^{h^{-1}}) \xi_n(\delta_T)(\Gamma g)>\nonumber\\
&=&\int_{\Gamma\backslash G} \int_K  \int_K  <\xi_n(\psi_T)(\Gamma
g),  \gamma(k^{-1}h)b(\Gamma g k^{-1}h)\gamma(k^{-1}h)^{-1} u(l^{h^{-1}}) \xi_n((\pi(k^{-1})\otimes \gamma(k^{-1}))\delta_T)(\Gamma
g)>\nonumber\\
&=&\int_{\Gamma\backslash G} \int_K  \int_K  <\xi_n(\psi_T)(\Gamma
g),  \gamma(k)^{-1}b(\Gamma gk)\gamma(k^{-1})u(l^{k^{-1}h}) \xi_n((\pi(h)\otimes \gamma(h))\delta_T)(\Gamma
g)>\label{eq24}
\end{eqnarray}

The Proposition \ref{pp8} now follows from the comparison of
(\ref{eq23}) and (\ref{eq24}).\hB

\subsubsection{}

\begin{lem}\label{ssk}
We have $\supp(\sigma_{micro})\subset q(\Gamma\backslash G\times \{l\})$.
\end{lem}
\proof
Let $f\in C(SM,p^*\End(V(\gamma))$ be such that $\tilde f_l=0$.
We must show that
$\sigma_{\xi_n(\psi_T)}(f)=o(1)$.

Let $D(G,\gamma)$ denote the algebra of $G$-invariant differential
operators on $G\times_KV_\gamma$. Note that the operators of $D(G,\gamma)$ descent to
$\Gamma\backslash G\times_KV_\gamma$.
The right-regular representation induces a homomorphism
$\tau:U(\gaaa)^K\rightarrow D(G,\gamma)$
such that
$(\tau(U)f)(g):=f(gU)$.
 If we compose $\tau$ with the symmetrization map (see \ref{lj1}),
then we get a linear map
$$D:Sym(\paaa)^K\rightarrow D(G,\gamma)\ .$$
Let $p\in Sym^d(\paaa)^K$. Then the symbol of the corresponding degree-$d$ differential operator
is given by the function
$s(D(p))\in C(SM,p^*\End(V_\gamma))$, $s(D(p))([\Gamma g,\lambda])=p(\lambda)$.
There is a Harish-Chandra homomorphism
$\Phi_\gamma:D(G,\gamma)\rightarrow U(\aaaa)\otimes \End_M(V_\gamma)$
(see \cite{olbrich}). We identify $U(\aaaa)\cong Sym(\aaaa)$ naturally.
For $D\in D(G,\gamma)$ we have (\cite{olbrich}, Lemma 2.13)
$$D\xi_n(\psi_T)=\xi_n(\psi_{\Phi_\gamma(D)(\lambda_n)\circ
  T})\ .$$
In addition, the $\End_M(V_\gamma)$-valued polynomial
$\aaaa^*_\C\ni\lambda\mapsto p(\lambda)\id-\Phi_\gamma(D)(\lambda)$ has degree
at most $d-1$ (\cite{olbrich}, Lemma 2.6).
It follows that
$$D(p)\xi_n(\psi_T)=p(\lambda_n)\xi_n(\psi_T)+ \|\lambda_n\|^do(1)\
.$$

Let now $f\in C(SM,p^*\End(V(\gamma))$.
Then we have
$$p(l)\sigma_{\xi_n(\psi_T)}(f)=\sigma_{\xi_n(\psi_T)}(s(D(p))f)+o(1)\ .$$ 
If $p(l)=0$, then $\sigma_{\xi_n(\psi_T)}(s(D(p))f)=o(1)$.
We now argue as in \cite{sv04}. If $\tilde f_l=0$, then it can be
approximated by products of the form
$s(D(p)) h$ with $p(l)=0$ (here one has to use the regularity of $l$ again). This implies the lemma. \hB

\subsubsection{}

Sending $[\Gamma g,l]$ to $\Gamma
gM$ identifies $q(\Gamma\backslash G\times\{l\})$ with the double quotient
$\Gamma\backslash G/M$.
In order to finish the proof of Theorem \ref{hgf} it remains to verify
the following lemma.

\begin{lem}\label{denss}
The symbols $s(D(d,U,b))$ with $d\ge 0$, $b\in C^K(\Gamma\backslash
G)\otimes \End(V_\gamma)$, and $U\in \cU^{\le d}(\gaaa)$,  span a dense subspace of
$C(q(\Gamma\backslash G\times\{l\}),p^*\End(V(\gamma)))\cong
C(\Gamma\backslash G\times_M \End(V_\gamma))$. \end{lem}
\proof
We have a $K$-bundle
$$\Gamma\backslash G\times K/M\times \End(V_\gamma)\rightarrow
\Gamma\backslash G\times_M\End(V_\gamma)$$
given by $(\Gamma g,kM,\Phi)\mapsto [\Gamma gk,\gamma(k^{-1})\Phi\gamma(k)]$. 
Since $l$ is regular, the functions
$k\mapsto u(l^{k^{-1}})$, $u\in S(\paaa)$,
span a dense subspace of $C(K/M)$ (see the proof of Lemma \ref{45d}).
Therefore the functions
$(\Gamma g,kM)\mapsto b(\Gamma g) u(l^{k^{-1}})$ span a dense subspace
of $C(\Gamma \backslash G\times K/M)\otimes \End(V_\gamma)$.
It follows that the $K$-averages 
$\Gamma g \mapsto \int_K \gamma(k) b(\Gamma g k)\gamma(k)^{-1} u(l^{k^{-1}})$
span a dense subspace of
$C(\Gamma\backslash G\times_M \End(V_\gamma))$. 
\hB

 \end{document}